\documentclass{amsart}
\usepackage{amsmath,amssymb}
\usepackage[dvips]{graphics}
\newtheorem{Theorem}{Theorem}[section]
\newtheorem{Lemma}[Theorem]{Lemma}

\newtheorem{Example}[Theorem]{Example}

\newtheorem{Remark}[Theorem]{Remark}
\newtheorem{Conjecture}[Theorem]{Conjecture}
\newtheorem{Question}[Theorem]{Question}

\makeatletter
\@addtoreset{figure}{section}%{subsection}
\def\@thmcountersep{-}
\makeatother

\numberwithin{equation}{section}

\begin{document} 

\title[Delta edge-homotopy invariants of spatial graphs]{Delta edge-homotopy invariants of spatial graphs via disk-summing the constituent knots}

%    Information for author
\author{Ryo Nikkuni}
\address{Department of Mathematics, Faculty of Education, Kanazawa University, Kakuma-machi, Kanazawa, Ishikawa, 920-1192, Japan}
\email{nick@ed.kanazawa-u.ac.jp}

%    General info

\date{}

\dedicatory{}


\begin{abstract}
In this paper we construct some invariants of spatial graphs 
by disk-summing the constituent knots and show the delta edge-homotopy 
invariance of them. As an application, 
we show that there exist infinitely many slice spatial embeddings of 
a planar graph up to delta edge-homotopy, and there exist infinitely many
boundary spatial embeddings of a planar graph up to delta edge-homotopy. 
\end{abstract}

\maketitle

\section{Introduction} 

Throughout this paper we work in the piecewise linear category. 
Let $G$ be a finite graph. 
An embedding of $G$ into the $3$-sphere is called a 
{\it spatial embedding} of $G$ or 
simply a {\it spatial graph}. A graph $G$ is said to be {\it planar} if 
there exists an embedding of $G$ into the $2$-sphere, and a 
spatial embedding of a planar graph $G$ is said to be {\it trivial} 
if it is ambient isotopic to 
an embedding of $G$ into the $2$-sphere in the $3$-sphere. 
Note that a trivial spatial 
embedding of a planar graph is unique up to ambient isotopy \cite{mason69}. 

\begin{figure}[htbp]
      \begin{center}
\scalebox{0.42}{\includegraphics*{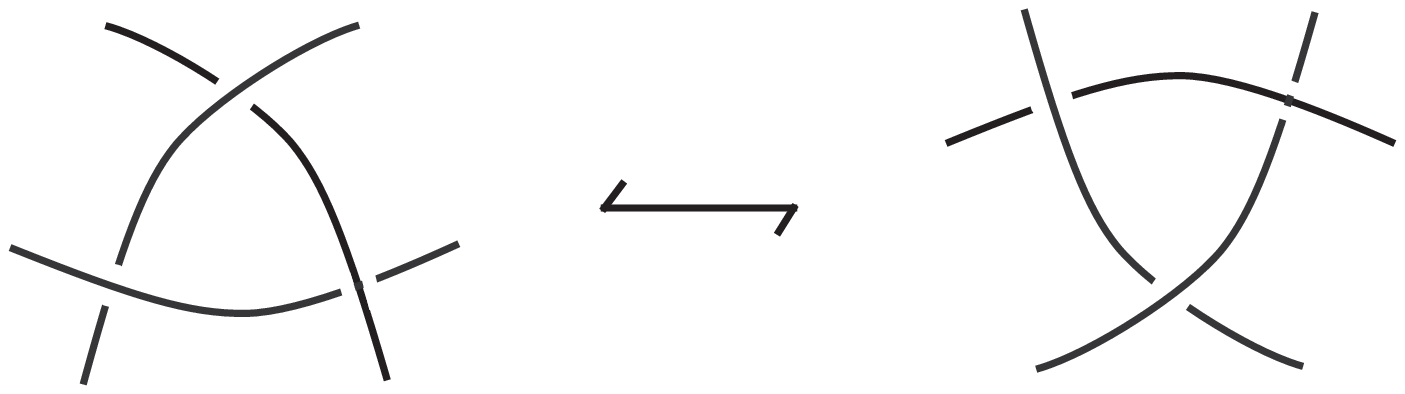}} 
      \end{center}
   \caption{}
  \label{delta}
\end{figure}

A {\it delta move} is a local deformation on a spatial graph as illustrated 
in Fig. \ref{delta} which is known as an unknotting operation 
\cite{matveev87}, \cite{m-n89}. 
A delta move is called a {\it self delta move} if all three strings in the 
move belong to the same spatial edge. 
Two spatial embeddings of a graph are said to be 
{\it delta edge-homotopic} 
if they are transformed into each other by self delta moves 
and ambient isotopies \cite{nikkuni02}. If the graph is homeomorphic to 
the disjoint union of $1$-spheres, then this equivalence relation 
coincides with  
{\it self $\Delta$-equivalence} \cite{shibuya96} 
(or {\it delta link homotopy} \cite{nakanishi02}) 
on oriented links. 

For self $\Delta$-equivalence on oriented links, Shibuya proposed 
the following conjectures in \cite{shibuya96} and \cite{shibuya00}. 

\begin{Conjecture}\label{shibuya_conj1}
{\rm (\cite{shibuya96})}
Two cobordant oriented links are self $\Delta$-equivalent. 
\end{Conjecture}
\begin{Conjecture}\label{shibuya_conj2}
{\rm (\cite{shibuya00})}
Any boundary link is self $\Delta$-equivalent to the trivial link.  
\end{Conjecture}

He gave the partially affirmative answers to the conjectures above 
at the same time. 
He showed that any {\it ribbon link} is self 
$\Delta$-equivalent to the trivial link \cite{shibuya96}, 
and any $2$-component boundary link is 
self $\Delta$-equivalent to the trivial link 
\cite[Theorem 4.6]{shibuya00}. 
But Nakanishi-Shibuya showed that there exists a  
$2$-component link such that 
it is not self $\Delta$-equivalent but cobordant 
to the Hopf link \cite[Claim 4.5]{n-s00}, 
namely they gave a negative answer to Conjecture \ref{shibuya_conj1}. 
Moreover, Nakanishi-Shibuya-Yasuhara showed that 
there exists a $3$-component link such that 
it is not self $\Delta$-equivalent but cobordant 
to the Borromean rings \cite[Proposition 1]{nsy}. 
Note that both the Hopf link and the Borromean rings are not slice. 
On the other hand, Conjecture \ref{shibuya_conj2} was solved affirmatively 
by Shibuya-Yasuhara \cite{shibuya-yasuhara07}. 

On the outcome of the results above, we investigate a more general case.   
A spatial embedding of a planar graph is said to be {\it slice} if 
it is cobordant\footnote{See \cite{taniyama94} for the precise definition 
of spatial graph-cobordism.} to the trivial spatial embedding. 
A spatial embedding of a graph is called a 
{\it $\partial$-spatial embedding} 
if all knots in the embedding bound Seifert surfaces simultaneously 
such that the interiors of the surfaces are mutually disjoint and disjoint 
from the image of the embedding \cite{n-s04}. 
If the graph is homeomorphic to the disjoint union of $1$-spheres, 
then this definition coincides with the definition of the boundary link. 
We note that any non-planar 
graph does not have a $\partial$-spatial embedding 
\cite[Corollary 1.3]{n-s04}. 
Then we ask the following questions. 

\begin{Question}\label{Q}
{\rm (1)} Is any slice spatial embedding of a planar graph 
delta edge-homotopic to the trivial spatial embedding? \\
{\rm (2)} Is any $\partial$-spatial embedding of a graph 
delta edge-homotopic to the trivial spatial embedding?
\end{Question}

In fact, for {\it spatial theta curves}, the affirmative answers to 
Question \ref{Q} (1) and (2) have already given 
by the author 
\cite[Corollary 1.3, Corollary 1.5]{nikkuni05}. But 
our purpose in this paper is to give the negative answers 
to the Questions \ref{Q} (1) and (2) as follows. 

\begin{Theorem}\label{A}
{\rm (1)} There exist infinitely many slice spatial embeddings of a graph 
up to delta edge-homotopy. \\
{\rm (2)} There exist infinitely many $\partial$-spatial embeddings of a graph 
up to delta edge-homotopy. 
\end{Theorem}

To accomplish this, 
we construct some invariants of spatial graphs 
by considering a disk-summing operation among the constituent knots in a 
spatial graph in section $2$, and 
show the delta edge-homotopy 
invariance of them in section $3$ 
(Theorems \ref{main_inv1} and \ref{main_inv2}). 
In section $4$, we give some remarkable examples which imply 
Theorem \ref{A}. Any of those examples is demonstrated by 
a spatial handcuff graph (see the next section) all of whose constituent links 
are trivial up to self $\Delta$-equivalence. Therefore our examples also imply that delta edge-homotopy on spatial graphs behaves quite differently than self $\Delta$-equivalence on links.

\begin{Remark}\label{mainrem}
{\rm (1) Recently, Question \ref{Q} (1) for oriented links was solved affirmatively by Yasuhara \cite[Corollary 1.9]{yasuhara07}. \\
(2) A {\it sharp move} 
is a local deformation on a spatial {\it oriented} graph as illustrated in 
Fig. \ref{sharp} 
which is also known as an unknotting operation \cite{murakami85}. 
A sharp move is called a {\it self sharp move} if all four strings in the 
move belong to the same spatial edge.  
Two spatial embeddings of a graph are said to be 
{\it sharp edge-homotopic} 
(or {\it self sharp-equivalent} \cite{n-s04}) 
if they are transformed into each other by self sharp moves 
and ambient isotopies \cite{nikkuni05b}.\footnote{
This equivalence relation does not depend on the edge 
orientations.}
It is known that two delta edge-homotopic spatial embeddings of a graph 
are sharp edge-homotopic  
\cite[Lemma 2.1 (2)]{nikkuni05b}. 
The author showed that two cobordant spatial embeddings of a graph are  
sharp edge-homotopic \cite[Lemma 2.2]{nikkuni05b}, and the author and 
R. Shinjo showed that any $\partial$-spatial embedding of a graph 
is sharp edge-homotopic  
to the trivial spatial embedding \cite[\sc Theorem 1.5 (1)]{n-s04}. 

\begin{figure}[htbp]
      \begin{center}
\scalebox{0.42}{\includegraphics*{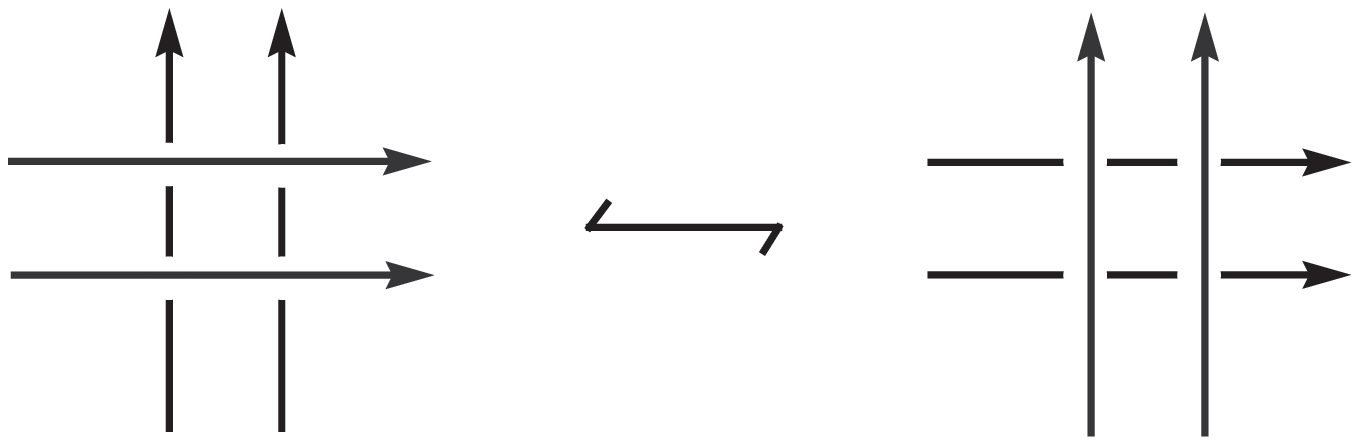}} 
   \end{center}
   \caption{}
  \label{sharp}
\end{figure} 
}
\end{Remark}
\section{Invariants} 

\begin{figure}[htbp]
      \begin{center}
\scalebox{0.32}{\includegraphics*{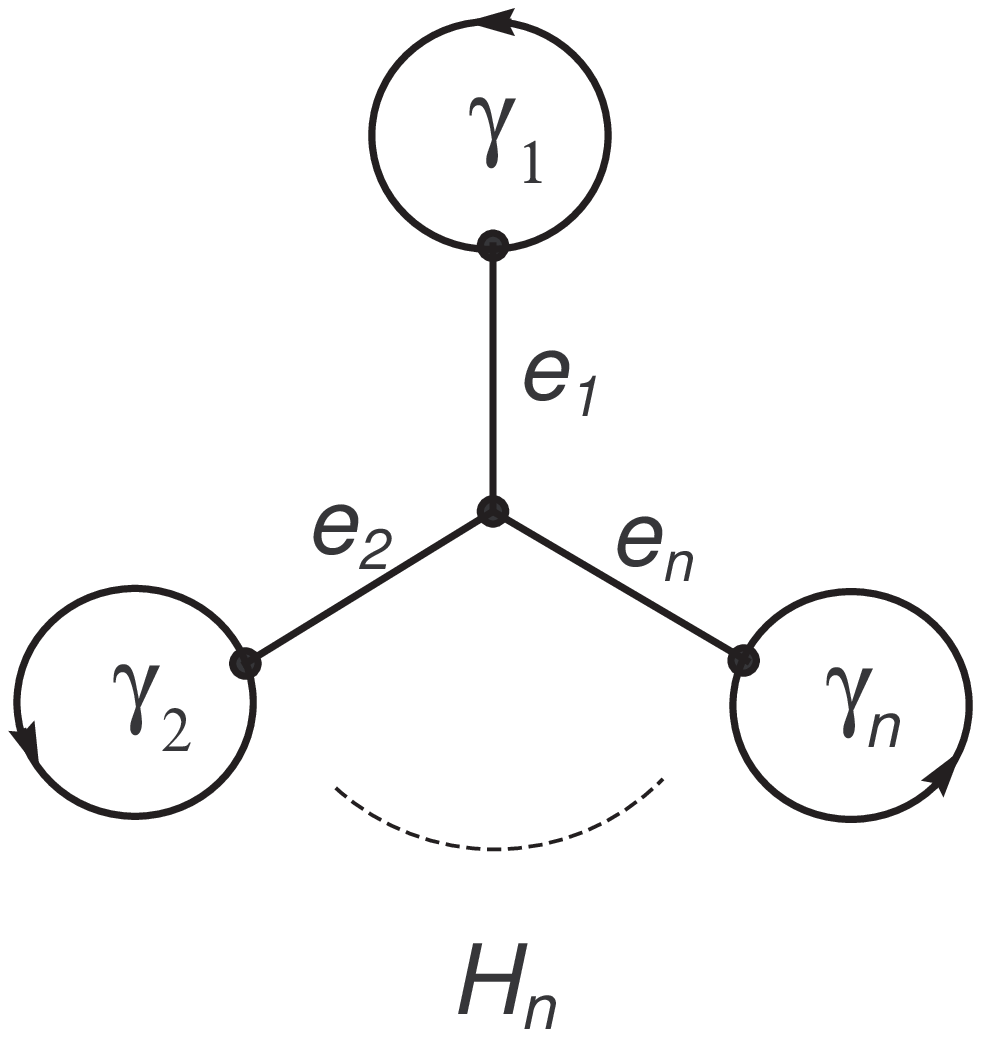}} 
      \end{center}
   \caption{}
  \label{3-handcuff}
\end{figure}

In this section we introduce the invariants of spatial 
graphs needed later. Let $H_{n}~(n\ge 2)$ be the graph as illustrated in 
Fig. \ref{3-handcuff}. We give the label to each of the edges 
and give an orientation to each of the loops as presented in 
Fig. \ref{3-handcuff}. 
A spatial embedding of $H_{n}$ 
is called a {\it spatial $n$-handcuff graph}, or simply a 
{\it spatial handcuff graph} if $n=2$. On that occasion, we regard 
$e_{1}\cup e_{2}$ as an edge of $H_{2}$ and denote by $e$. 

Let $L=J_{1}\cup J_{2}\cup\cdots\cup J_{n}$ be an ordered and oriented 
$n$-component link. Let $D$ be an oriented $2$-disk and  
$x_{1},x_{2},\ldots,x_{n}$ are mutually 
disjoint arcs in $\partial D$, where $\partial D$ has the orientation 
induced by the one of $D$, and 
these arcs 
appear along the orientation of $\partial D$ in order and 
each arc has an orientation induced 
by the one of $\partial D$. 
We assume that $D$ is embedded in the $3$-sphere so that 
$D\cap L=x_{1}\cup x_{2}\cup\cdots \cup x_{n}$ and $x_{i}\subset J_{i}$ with 
opposite orientations for any $i$. Then we call a knot 
$K_{D}^{12\cdots n}=L\cup \partial D
-\cup_{i=1}^{n}{\rm int}x_{i}$ a 
{\it $D$-sum} of $L$. For a spatial $n$-handcuff graph $f$, 
we denote $f(\gamma_{1}\cup \gamma_{2}\cup \cdots\cup \gamma_{n})$ 
by $L_{f}$ and 
consider a $D$-sum of $L_{f}$ so that 
$f(e_{1}\cup e_{2}\cup \cdots\cup e_{n})\subset D$ and 
$f(e_{i})\cap \partial D=f(e_{i}\cap \gamma_{i})
\subset {\rm int}x_{i}$ for any $i$. 
We call such a $D$-sum of $L_{f}$ a {\it $D$-sum of $L_{f}$ with respect 
to $f$} and denote it by $K_{D}^{12\cdots n}(f)$. 

For a spatial handcuff graph $f$, we define that 
\begin{eqnarray*}
{n}_{12}(f,D)=
a_{2}(K_{D}^{12}(f))-a_{2}(f(\gamma_{1}))-a_{2}(f(\gamma_{2}))
\end{eqnarray*}
and denote the modulo ${\rm lk}(L_{f})$ reduction of ${n}_{12}(f,D)$ by 
$\bar{n}_{12}(f)$, where ${\rm lk}$ denotes the {\it linking number} in the 
$3$-sphere. Then we have the following. 

\begin{Theorem}\label{main_inv1}
If two spatial handcuff graphs $f$ and $g$ 
are delta edge-homotopic, 
then $\bar{n}_{12}(f)=\bar{n}_{12}(g)$. 
\end{Theorem}

On the other hand, let $f$ be a spatial 3-handcuff graph and 
$K_{D}^{123}(f)$ a $D$-sum of $L_{f}$ with respect to $f$. Then, 
by using the same disk $D$, we can obtain 
three knots $K_{D}^{12}(f)$, $K_{D}^{23}(f)$ and 
$K_{D}^{13}(f)$ by forgetting the components $f(\gamma_{3})$, 
$f(\gamma_{1})$ and $f(\gamma_{2})$, respectively, namely  
by the $D$-sums of sublinks 
$f(\gamma_{1})\cup f(\gamma_{2})$, $f(\gamma_{2})\cup f(\gamma_{3})$ and 
$f(\gamma_{1})\cup f(\gamma_{3})$ of $L_{f}$. 
Then we define that 
\begin{eqnarray*}
{n}_{123}(f,D)=
-v_{3}(K_{D}^{123}(f))
+\sum_{1\le i<j\le 3}v_{3}(K_{D}^{ij}(f))
-\sum_{i=1}^{3}v_{3}(f(\gamma_{i})), 
\end{eqnarray*}
where 
$v_{3}(J)=(1/36)V_{J}^{(3)}(1)$ and $V_{J}^{(3)}(1)$ denotes the third 
derivative at $1$ of the {\it Jones polynomial} \footnote{
We calculate the Jones polynomial of a knot by the skein relation 
\begin{eqnarray*}
tV_{J_{+}}(t)-t^{-1}V_{J_{-}}(t)
=(t^{-\frac{1}{2}}-t^{\frac{1}{2}})V_{J_{0}}(t),
\end{eqnarray*} 
where $J_{+}$ and $J_{-}$ are two oriented knots and $J_{0}$ an 
oriented $2$-component link which are identical except inside 
the depicted regions  
as illustrated in Fig. \ref{skein}.} of a knot $J$. 
Assume that $L_{f}$ is {\it algebraically split}, namely all of the pairwise 
linking numbers of $L_{f}$ are zero. 
Then we denote the modulo $\mu_{123}(L_{f})$ reduction of ${n}_{123}(f,D)$ 
by $\bar{n}_{123}(f)$, where 
$\mu_{123}$ denotes the {\it triple linking number}, namely {\it Milnor's $\mu$-invariant} of length $3$ of a 
$3$-component algebraically split link \cite{milnor54}, \cite{milnor57}. 
Then we have the following. 

\begin{Theorem}\label{main_inv2}
Let $f$ and $g$ be two spatial $3$-handcuff graphs which are delta edge-homotopic. Assume that both $L_{f}$ and $L_{g}$ are algebraically split. Then it holds that $\bar{n}_{123}(f)=\bar{n}_{123}(g)$. 
\end{Theorem}

For example, if a spatial handcuff (resp. $3$-handcuff) graph $f$ contains 
a Hopf link (resp. Borromean rings), then our invariants are no use. 
But if $L_{f}$ is {\it link-homotopic} \cite{milnor54} to the trivial link, 
then our invariants 
take effect on our purpose. Because 
any slice link is link-homotopic to the trivial link \cite{giffen79}, 
\cite{goldsmith79}, 
and any boundary link is also link-homotopic to the trivial link 
\cite{c-f88}, \cite{dimovski88}. 
We prove Theorems \ref{main_inv1} and \ref{main_inv2} in the next section. 

\section{Proofs of Theorems \ref{main_inv1} and Theorems \ref{main_inv2}} 

To prove Theorems \ref{main_inv1} and \ref{main_inv2}, we first recall 
some results and show a lemma needed later. 

\begin{figure}[htbp]
      \begin{center}
\scalebox{0.35}{\includegraphics*{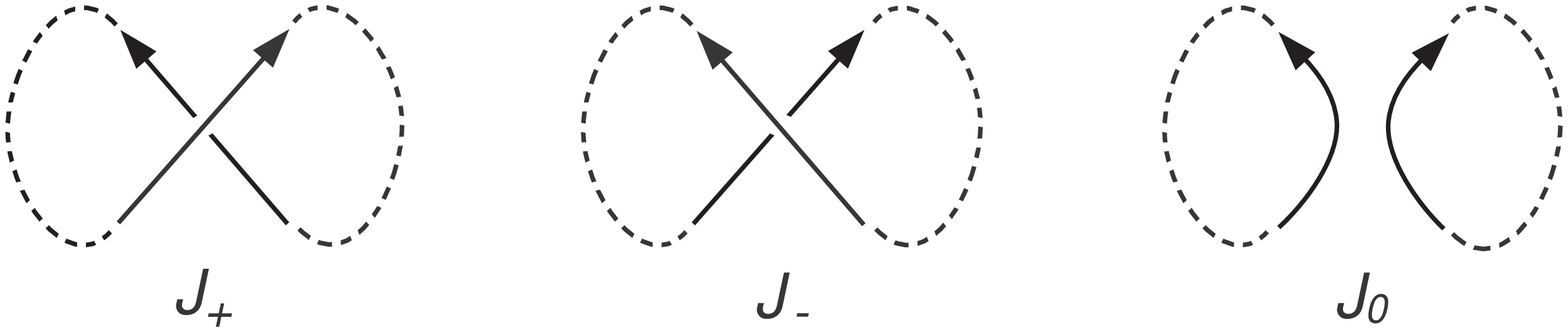}} 
      \end{center}
   \caption{}
  \label{skein}
\end{figure}

\begin{figure}[htbp]
      \begin{center}
\scalebox{0.35}{\includegraphics*{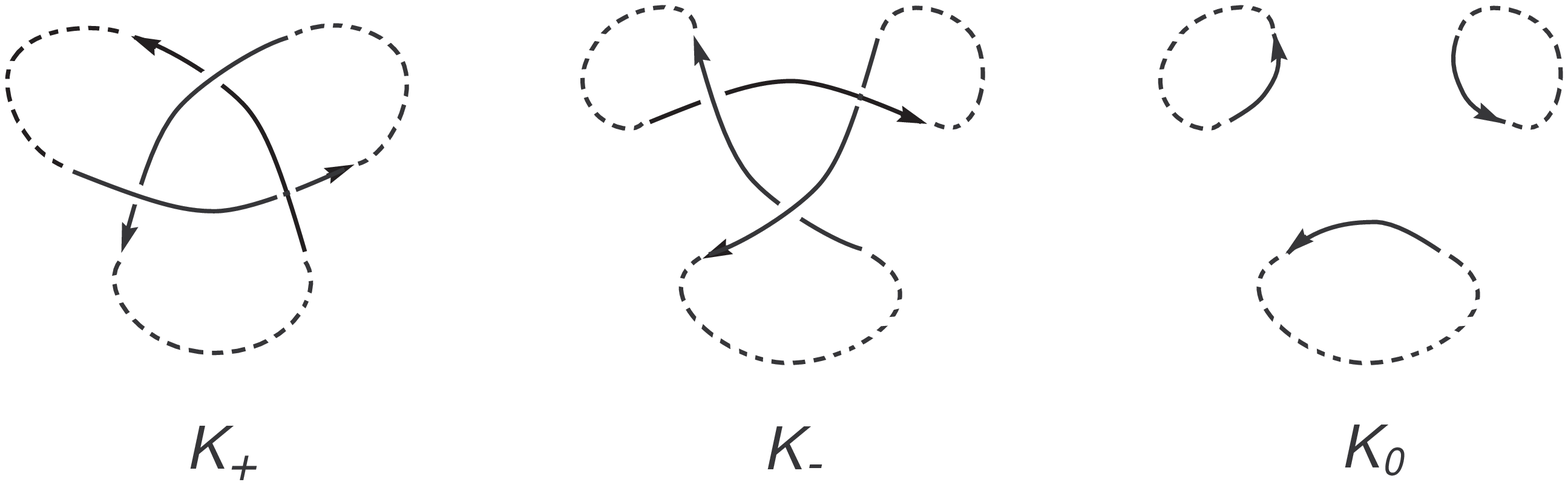}}
      \end{center}
   \caption{}
  \label{d_skein_triple}
\end{figure}

\begin{Lemma}\label{vlemma2}   
Let $J_{+}$ and $J_{-}$ be two oriented knots and $J_{0}=K_{1}\cup K_{2}$ an 
oriented $2$-component link which are identical except inside 
the depicted regions  
as illustrated in Fig. \ref{skein}. Then we have that 

\vspace{0.2cm}
{\rm (1) (\cite[Lemma~5.6]{kaufmann83})} 
$a_{2}(J_{+})-a_{2}(J_{-})={\rm lk}(J_{0})$.

\vspace{0.2cm}
{\rm (2) (\cite[Proposition~4.2]{nikkuni05b})} \\ 
$V_{J_{+}}^{(3)}(1)-V_{J_{-}}^{(3)}(1)
=36a_{2}(J_{+})+18\{{\rm lk}(J_{0})\}^{2}
-36\{a_{2}(K_{1})+a_{2}(K_{2})\}$. 
\end{Lemma}
\begin{Lemma}\label{vlemma3}   
Let $K_{+}$ and $K_{-}$ be two oriented knots and 
$K_{0}$ an 
oriented $3$-component link which are identical except inside 
the depicted regions 
as illustrated in Fig. \ref{d_skein_triple}. Then we have that 

\vspace{0.2cm}
{\rm (1) (\cite[Theorem 1.1]{okada90})} 
$a_{2}(K_{+})-a_{2}(K_{-})=1$.

\vspace{0.2cm}
{\rm (2) (\cite[Theorem 3.2]{nikkuni02})} 
$\displaystyle V_{K_{+}}^{(3)}(1)-V_{K_{-}}^{(3)}(1)
=36{\rm Lk}(K_{0})-18$, 
\vspace{0.2cm}\\
where ${\rm Lk}(L)$ denotes the {\it total linking number} of 
an oriented link 
$L$. 
\end{Lemma}
\begin{Lemma}\label{delta_lem}
{\rm (1)} Let $f$ be a spatial handcuff graph. Then any of the 
self delta moves 
on $f(e)$ is realized by self delta moves on $f(\gamma_{1})$ and 
ambient isotopies. \\
{\rm (2)} Let $f$ be a spatial $3$-handcuff graph. Then any of the 
self delta moves  
on $f(e_{i})$ is realized by self delta moves on $f(\gamma_{i})$ and 
ambient isotopies $(i=1,2,3)$.
\end{Lemma}

\begin{figure}[htbp]
      \begin{center}
    \scalebox{0.32}{\includegraphics*{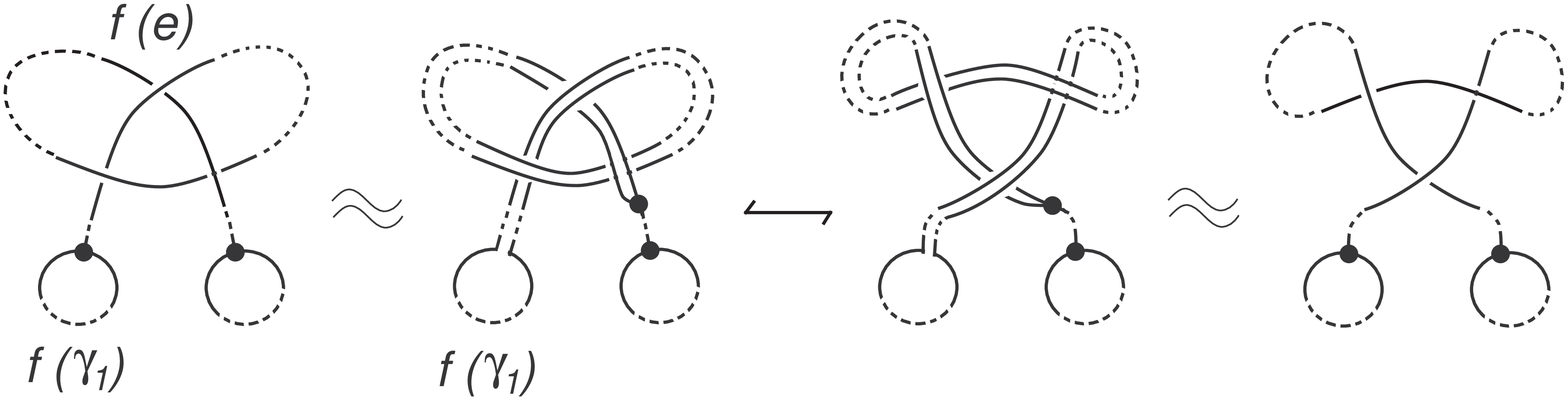}}
      \end{center}
   \caption{}
  \label{doubled_delta}
\end{figure}

{\it Proof.} 
(1) We can see that a self delta move on $f(e)$ is realized 
by a ``doubled-delta move'' on $f(\gamma_{1})$, see Fig. \ref{doubled_delta}. 
It is easy to see that a doubled-delta move is realized by eight delta 
moves on the strings in the move and ambient isotopies. Thus we have the 
result. \\
(2) We can show in the same way as (1). \hfill $\square$

\begin{figure}[htbp]
      \begin{center}
\scalebox{0.32}{\includegraphics*{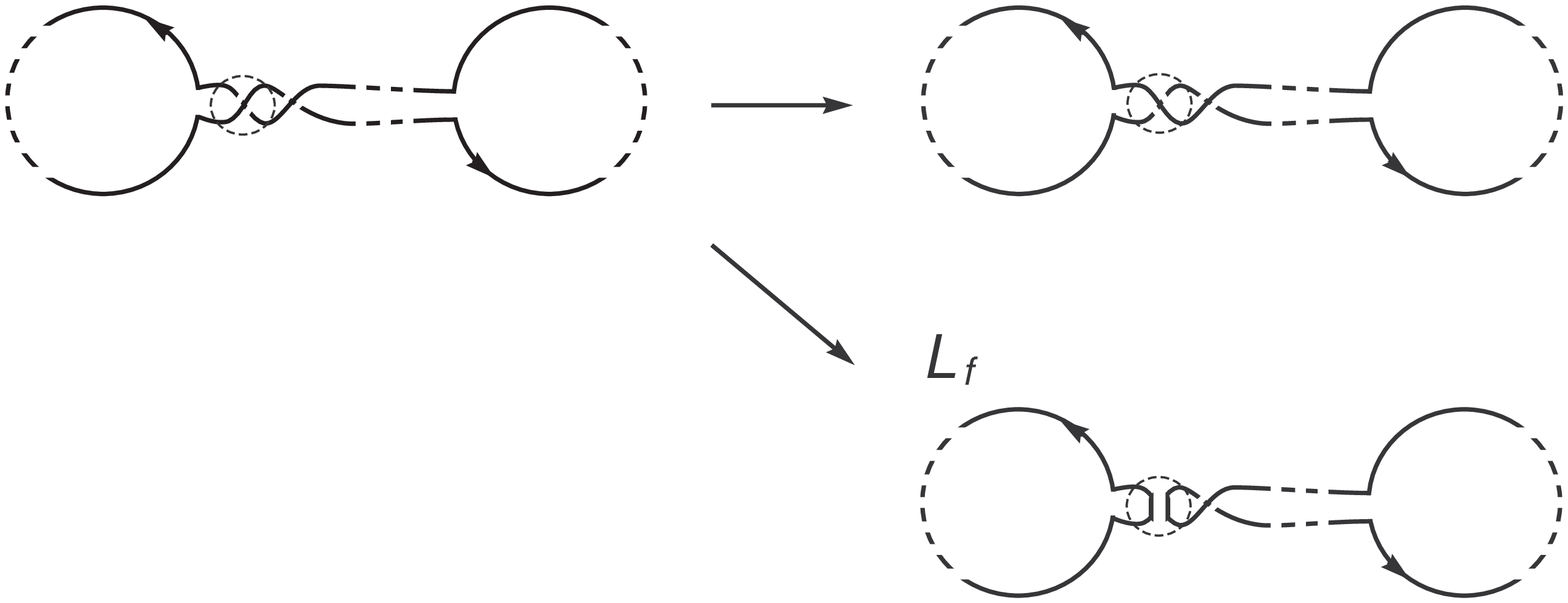}}
      \end{center}
   \caption{}
  \label{skein_inv}
\end{figure}

\vspace{0.2cm}
{\it Proof of Theorem \ref{main_inv1}.}  
We first show that $\bar{n}_{12}(f)$ is an ambient isotopy invariant. 
Let $K_{D}^{12}(f)$ be a $D$-sum of $L_{f}$ with respect to $f$ 
and $K_{D'}^{12}(f)$ another $D'$-sum 
of $L_{f}$ with respect to $f$. 
We may assume that $K_{D'}^{12}(f)$ is obtained from 
$K_{D}^{12}(f)$ 
by a positive full twist of the band corresponding to $f(e)$. 
Then by Lemma \ref{vlemma2} (1) we have that 
\begin{eqnarray*}
n_{12}(f,D')-n_{12}(f,D)&=&a_{2}(K_{D'}^{12}(f))-a_{2}(K_{D}^{12}(f))
={\rm lk}(L_{f}), 
\end{eqnarray*}
see Fig. \ref{skein_inv}. This implies that 
$\bar{n}_{12}(f)$ is an ambient isotopy invariant. 

\begin{figure}[htbp]
      \begin{center}
\scalebox{0.4}{\includegraphics*{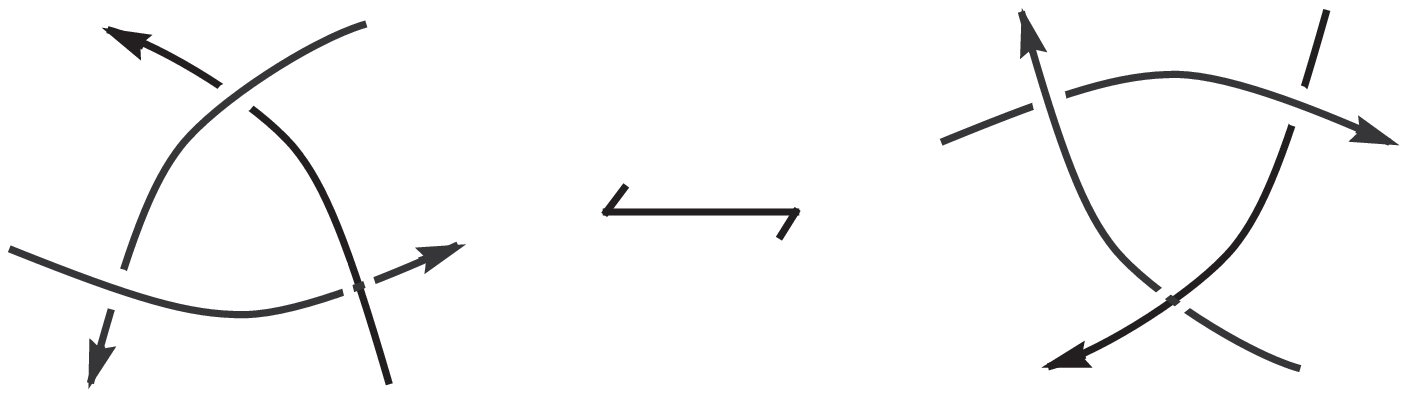}}      
      \end{center}
   \caption{}
  \label{oriented_delta}
\end{figure}

\begin{figure}[htbp]
      \begin{center}
\scalebox{0.32}{\includegraphics*{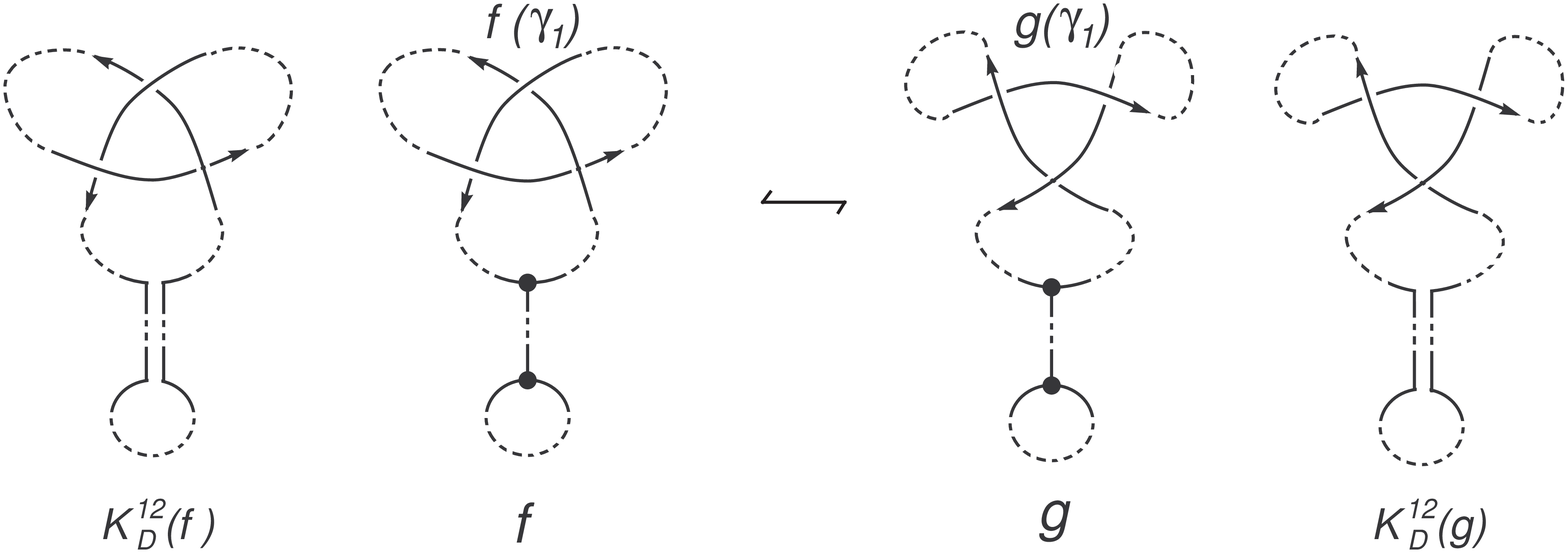}} 
      \end{center}
   \caption{}
  \label{selfdelta_inv1}
\end{figure}

Next we show that $\bar{n}_{12}(f)$ is a delta edge-homotopy invariant. 
Let $f$ and $g$ be delta edge-homotopic two spatial handcuff graphs.  
Then, by Lemma \ref{delta_lem} (1), $g$ is obtained from $f$ 
by self delta moves on $f(\gamma_{i})$ ($i=1,2$) and ambient isotopies. 
Moreover it is known that each of the oriented delta moves can be realized by 
the one as illustrated in 
Fig. \ref{oriented_delta} \cite[Fig. 1.1]{m-n89}.
Hence we may assume that $g$ is obtained from $f$ by a self delta move 
on $f(\gamma_{1})$ as illustrated in Fig. \ref{selfdelta_inv1}  
without loss of generality. 
Let $K_{D}^{12}(f)$ be the $D$-sum of $L_{f}$ with respect to $f$ 
as illustrated in 
Fig. \ref{selfdelta_inv1} and $K_{D}^{12}(g)$ the $D$-sum of 
$L_{g}$ with respect to $g$ by using the same $D$ 
as illustrated in Fig. \ref{selfdelta_inv1}. 
Namely $K_{D}^{12}(f)$ and $K_{D}^{12}(g)$ are identical except 
the depicted parts which represents the delta move. Note that 
$f(\gamma_{2})$ and $g(\gamma_{2})$ are ambient isotopic. 
Then by Lemma \ref{vlemma3} (1) we have that 
\begin{eqnarray*}
&&n_{12}(f,D)-n_{12}(g,D)\\
&=&a_{2}(K_{D}^{12}(f))-a_{2}(K_{D}^{12}(g))-
\left\{a_{2}(f(\gamma_{1}))-a_{2}(g(\gamma_{1}))\right\}\\
&=&1-1=0.
\end{eqnarray*}
Since a delta move preserves the linking number, we have that 
$\bar{n}_{12}(f)=\bar{n}_{12}(g)$. This completes the proof. \hfill $\square$

\begin{Remark}\label{rem_inv}   
{\rm (1) By the first half of the proof of 
Theorem \ref{main_inv1}, we can also see that the modulo 
${\rm lk}(L_{f})$ reduction of 
$a_{2}(K_{D}^{12}(f))$ is an ambient isotopy invariant of a 
spatial handcuff graph $f$.  \\
(2) For a spatial $n$-handcuff graph $f$ and a $D$-sum 
$K_{D}^{12\cdots n}(f)$ of $L_{f}$ with respect to $f$, 
we can generalize Theorem \ref{main_inv1} as follows. 
Let $l$ be the greatest common divisor of 
${\rm lk}(f(\gamma_{i}),L_{f}-f(\gamma_{i}))~(i=1,2,\ldots,n)$. 
Then it can be shown that the modulo $l$ reduction of 
\begin{eqnarray*}
a_{2}(K_{D}^{12\cdots n}(f))-\sum_{i=1}^{n}a_{2}(f(\gamma_{i}))
\end{eqnarray*}
is a delta edge-homotopy invariant of $f$ 
in the same way as the proof of Theorem 
\ref{main_inv1}. But this generalized 
version for $n\ge 3$ is not so strong as we will see in 
Examples \ref{ex3} and \ref{ex5}. 
}
\end{Remark}

To prove Theorem \ref{main_inv2}, we recall another result. 
By Polyak's formula of the triple linking number \cite{polyak97}, we have the following. 

\begin{Lemma}\label{polyak} 
{\rm (\cite{polyak97})} 
Let $L=J_{1}\cup J_{2}\cup J_{3}$ be an ordered and oriented algebraically split $3$-component link. Let $K_{D}$ be a $D$-sum of $L$ and 
$K_{D}^{23},K_{D}^{13}$ and $K_{D}^{12}$ three knots obtained from 
$K_{D}^{123}$ by forgetting the components 
$J_{1}, J_{2}$ and $J_{3}$, respectively. 
Then it holds that 
\begin{eqnarray*}
{\mu}_{123}(L)= -a_{2}(K_{D}^{123})
+\sum_{1\le i<j\le 3}a_{2}(K_{D}^{ij})
-\sum_{i=1}^{3}a_{2}(J_{i}).
\end{eqnarray*}  
\end{Lemma}

\begin{figure}[htbp]
      \begin{center}
\scalebox{0.35}{\includegraphics*{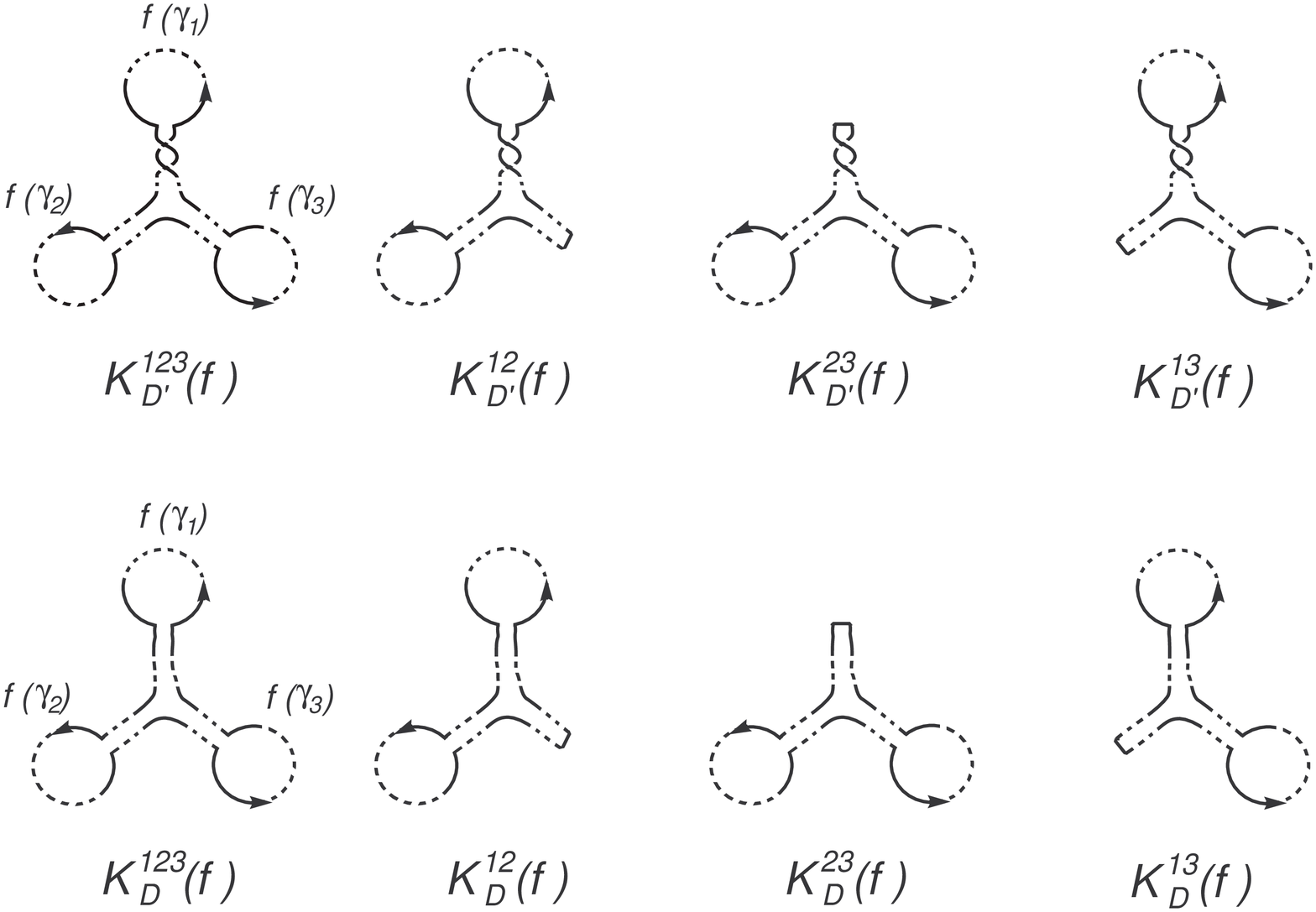}} 
      \end{center}
   \caption{}
  \label{hand_quad}
\end{figure}

\begin{figure}[htbp]
      \begin{center}
\scalebox{0.35}{\includegraphics*{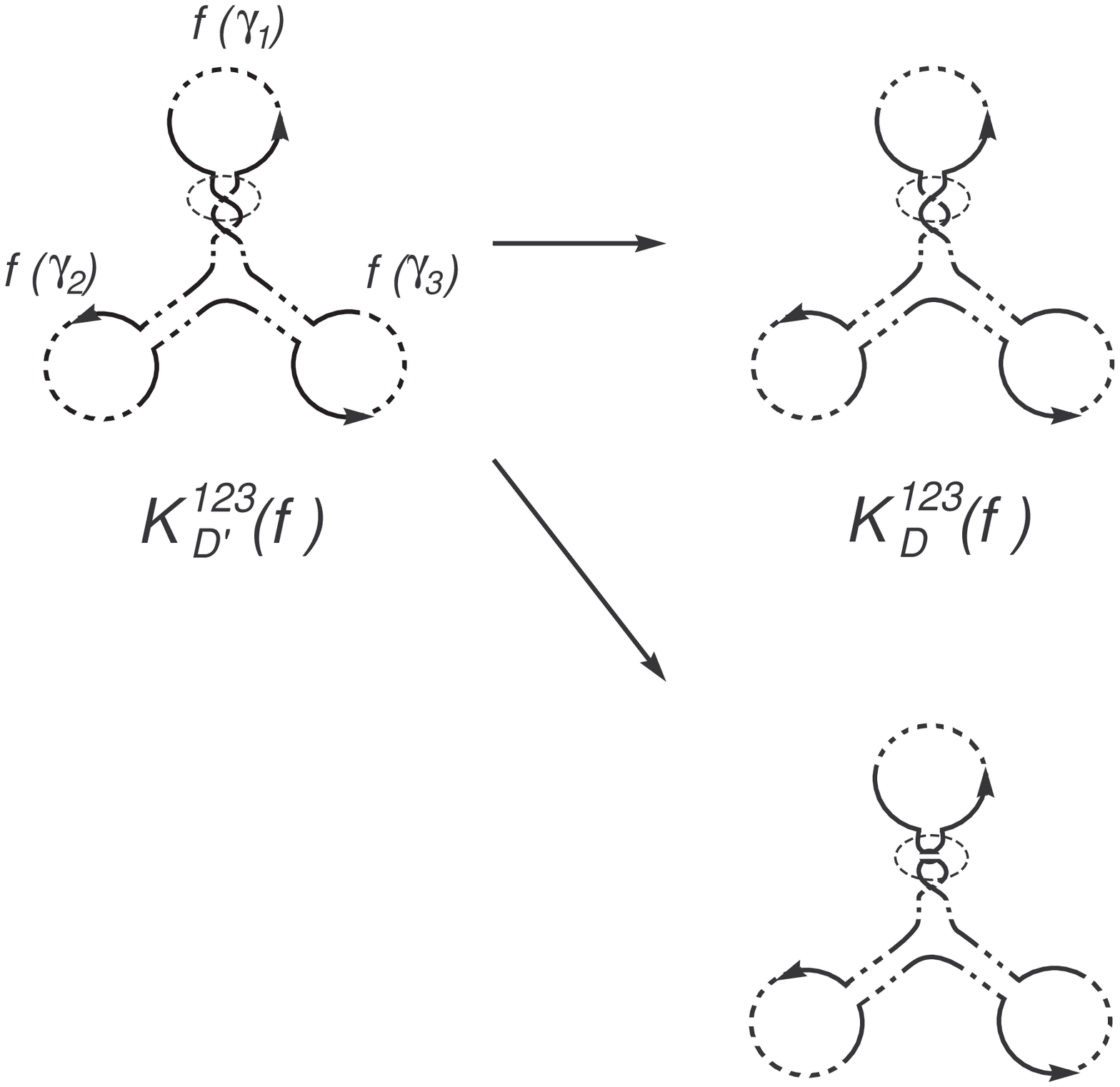}}
      \end{center}
   \caption{}
  \label{hand_quad2}
\end{figure}

{\it Proof of Theorem \ref{main_inv2}.} 
We first show that $\bar{n}_{123}(f)$ is an ambient isotopy invariant. 
Let $K_{D}^{123}(f)$ be a $D$-sum of $L_{f}$ with respect to $f$ 
and $K_{D'}^{123}(f)$ another $D'$-sum 
of $L_{f}$ with respect to $f$. 
We may assume that $K_{D'}^{123}(f)$ is obtained from 
$K_{D}^{123}(f)$ 
by a positive full twist of the band corresponding to $f(e_{1})$, 
see Fig. \ref{hand_quad}. Then by the skein relation as illustrated 
in Fig. \ref{hand_quad2}, 
Lemmas \ref{vlemma2} (2) and \ref{polyak}, 
we have that 
\begin{eqnarray*}
&&n_{123}(f,D')-n_{123}(f,D)\\
&=&-\left\{
v_{3}(K_{D'}^{123}(f))-v_{3}(K_{D}^{123}(f))
\right\}
+\left\{
v_{3}(K_{D'}^{12}(f))-v_{3}(K_{D}^{12}(f))
\right\}\\
&&+\left\{
v_{3}(K_{D'}^{13}(f))-v_{3}(K_{D}^{13}(f))
\right\}\\
&=&-a_{2}(K_{D'}^{123}(f))
-\frac{1}{2}{\rm lk}(f(\gamma_{1}),f(\gamma_{2})\cup f(\gamma_{3}))^{2}
+\left\{a_{2}(f(\gamma_{1}))+a_{2}(K_{D'}^{23}(f))\right\}\\
&&+a_{2}(K_{D'}^{12}(f))
+\frac{1}{2}{\rm lk}(f(\gamma_{1}),f(\gamma_{2}))^{2}
-\left\{a_{2}(f(\gamma_{1}))+a_{2}(f(\gamma_{2}))\right\}\\
&&+a_{2}(K_{D'}^{13}(f))
+\frac{1}{2}{\rm lk}(f(\gamma_{1}),f(\gamma_{3}))^{2}
-\left\{a_{2}(f(\gamma_{1}))+a_{2}(f(\gamma_{3}))\right\}\\
&=&-a_{2}(K_{D'}^{123}(f))+\sum_{1\le i<j \le 3}a_{2}(K_{D'}^{ij}(f))
-\sum_{i=1}^{3}a_{2}(f(\gamma_{i}))\\
&=& {\mu}_{123}(L_{f}). 
\end{eqnarray*}
Hence we have that $\bar{n}_{123}(f)$ is an ambient isotopy invariant. 


\begin{figure}[htbp]
      \begin{center}
\scalebox{0.32}{\includegraphics*{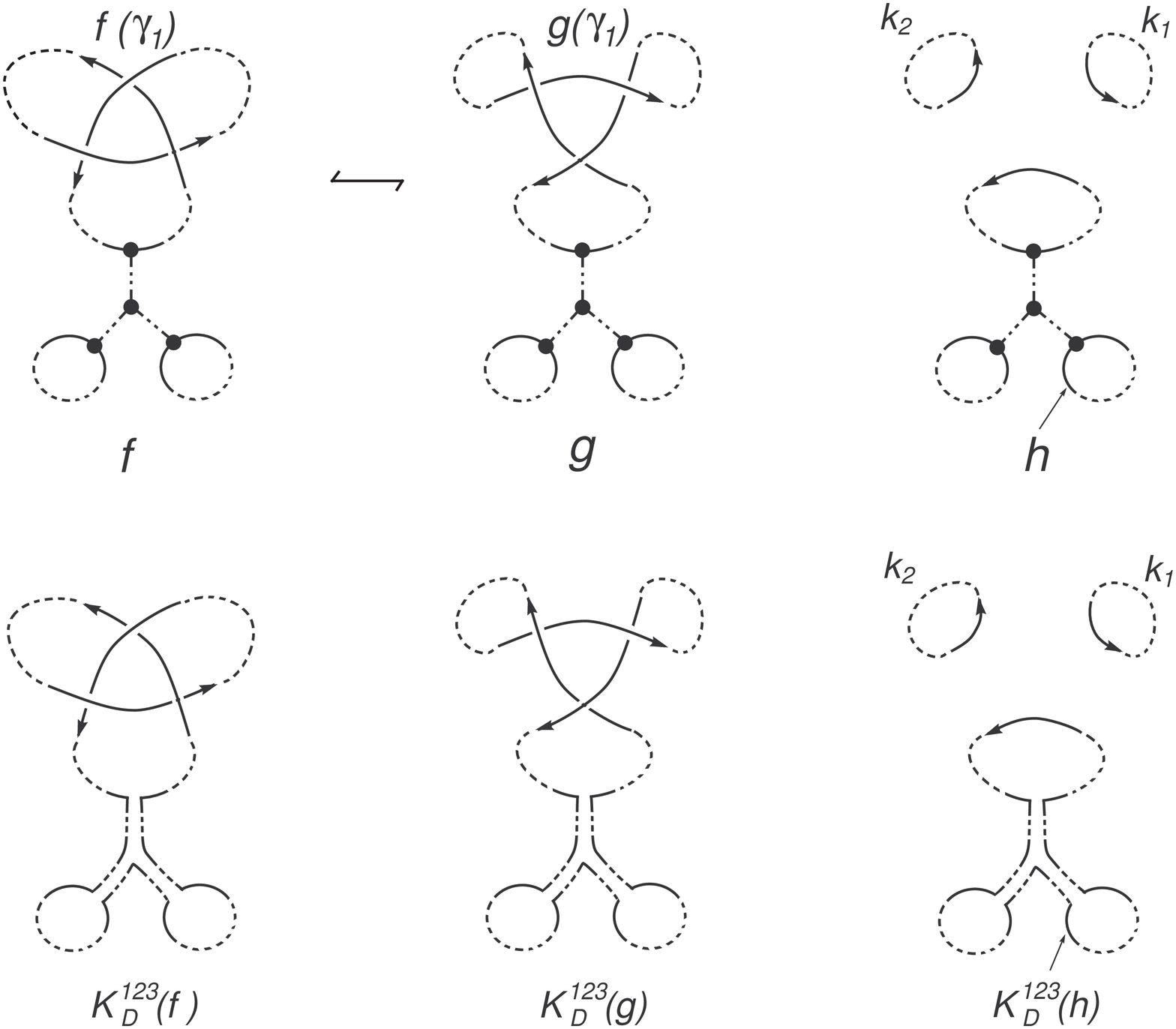}}
      \end{center}
   \caption{}
  \label{selfdelta_inv2}
\end{figure}

Next we show that $\bar{n}_{123}(f)$ is a delta edge-homotopy invariant. 
Let $f$ and $g$ be delta edge-homotopic two spatial $3$-handcuff graphs.  
Then, by Lemma \ref{delta_lem} (2), $g$ is obtained from $f$ 
by self delta moves on $f(\gamma_{i})$ ($i=1,2,3$) and ambient isotopies. 
Hence we may assume that $g$ is obtained from $f$ by a self delta move 
on $f(\gamma_{1})$ as illustrated in Fig. \ref{selfdelta_inv2}  
without loss of generality. 
Let $K_{D}^{123}(f)$ be the $D$-sum of $L_{f}$ with respect to $f$ 
as illustrated in 
Fig. \ref{selfdelta_inv2} and $K_{D}^{123}(g)$ the $D$-sum of 
$L_{g}$ with respect to $g$ 
by using the same $D$ as illustrated in Fig. \ref{selfdelta_inv2}. 
Namely $K_{D}^{123}(f)$ and $K_{D}^{123}(g)$ are identical except 
the depicted parts which represents the delta move. 
Let $h$ be the spatial $3$-handcuff graph and $k_{1}$ and $k_{2}$ two 
oriented knots as illustrated in Fig. \ref{selfdelta_inv2}, where 
$f(H_{3})$, $g(H_{3})$ and $h(H_{3})\cup k_{1}\cup k_{2}$ are identical 
except the depicted parts. 
Let $K_{D}^{123}(h)$ be the $D$-sum of $L_{h}$ with respect to $h$ 
by using the same $D$ as illustrated in Fig. \ref{selfdelta_inv2}. 
Then by Lemma \ref{vlemma3} (2) and the homological invariance of the 
linking number, we have that 
\begin{eqnarray*}
&&n_{123}(f,D)-n_{123}(g,D)\\
&=&
-{\rm lk}(k_{1},K_{D}^{123}(h))-{\rm lk}(k_{2},K_{D}^{123}(h))
-{\rm lk}(k_{1},k_{2})+\frac{1}{2}\\
&&+{\rm lk}(k_{1},K_{D}^{12}(h))+{\rm lk}(k_{2},K_{D}^{12}(h))
+{\rm lk}(k_{1},k_{2})-\frac{1}{2}\\
&&+{\rm lk}(k_{1},K_{D}^{13}(h))+{\rm lk}(k_{2},K_{D}^{13}(h))
+{\rm lk}(k_{1},k_{2})-\frac{1}{2}\\
&&-{\rm lk}(k_{1},h(\gamma_{1}))-{\rm lk}(k_{2},h(\gamma_{2}))
-{\rm lk}(k_{1},k_{2})+\frac{1}{2}\\
&=&0.
\end{eqnarray*}

Note that $L_{f}$ and $L_{g}$ are self $\Delta$-equivalent. Thus they are 
also link-homotopic, 
namely ${\mu}_{123}(L_{f})={\mu}_{123}(L_{g})$. 
Thus we have that $\bar{n}_{123}(f)=\bar{n}_{123}(g)$. 
This completes the proof. \hfill $\square$

\section{Examples} 
\begin{Example}\label{ex2}
{\rm Let $f_{m}$ be the spatial handcuff graph for $m\in {\mathbb N}$ as 
illustrated in Fig. \ref{slice_handcuff1}.
We can see that $L_{f_{m}}$ is the trivial $2$-component link for 
any $m\in {\mathbb N}$, namely ${\rm lk}(L_{f_{m}})=0$. 
We can also see that $f_{m}$ is slice 
by the hyperbolic transformation on $f_{m}(\gamma_{2})$ along the band $B$ shown in Fig. \ref{slice_handcuff1}. 

\begin{figure}[htbp]
      \begin{center}
\scalebox{0.375}{\includegraphics*{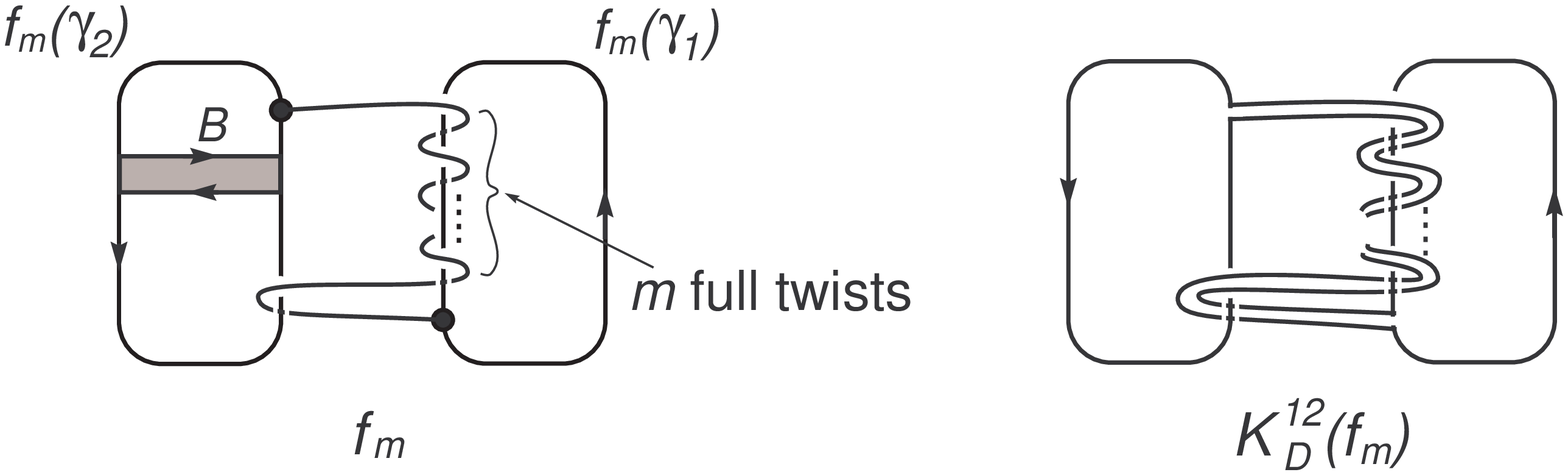}}
      \end{center}
   \caption{}
  \label{slice_handcuff1}
\end{figure} 

Now we consider the $D$-sum of $L_{f_{m}}$ with 
respect to $f_{m}$ as illustrated in 
Fig. \ref{slice_handcuff1}. 
Then by a calculation we have that $a_{2}(K_{D}^{12}(f_{m}))=2m$ and therefore 
$\bar{n}_{12}(f_{m})=2m$. Thus by Theorem \ref{main_inv1},  
we have that $f_{m}$ is not delta edge-homotopic to the trivial spatial 
handcuff graph for any $m\in {\mathbb N}$, and 
$f_{i}$ and $f_{j}$ are 
not delta edge-homotopic for $i\neq j$. 
}
\end{Example}
\begin{Example}\label{ex3}
{\rm Let $f_{m}$ be the spatial $3$-handcuff graph for $m\in {\mathbb N}$ as 
illustrated in Fig. \ref{slice_handcuff3}.
We can see that $L_{f_{m}}$ is the trivial $3$-component link for any $m\in {\mathbb N}$, namely $\mu_{123}(L_{f_{m}})=0$. We can also see that 
$f_{m}$ is slice by the hyperbolic transformation on $f_{m}(\gamma_{3})$ along the band $B$ shown in Fig. \ref{slice_handcuff3}. 

\begin{figure}[htbp]
      \begin{center}
\scalebox{0.35}{\includegraphics*{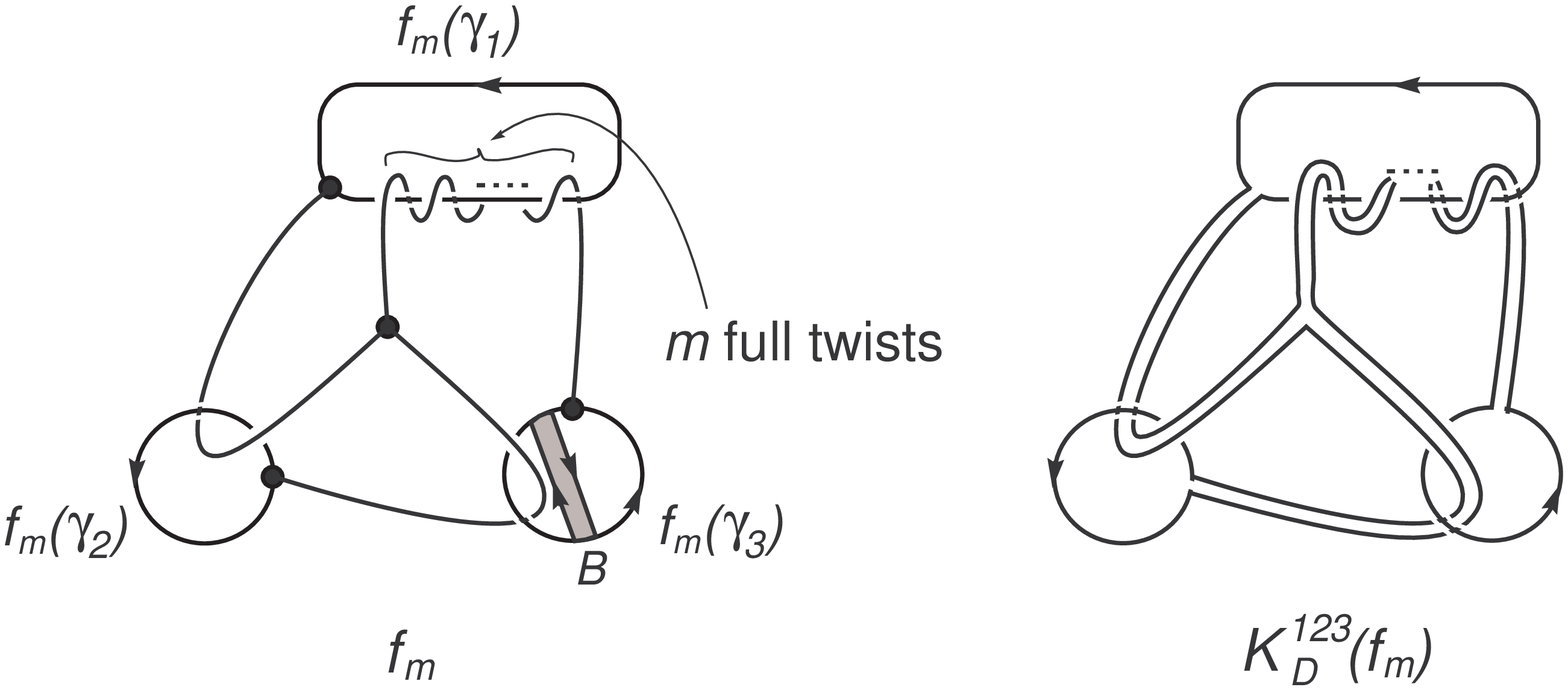}}
      \end{center}
   \caption{}
  \label{slice_handcuff3}
\end{figure} 

Now we consider the $D$-sum of $L_{f_{m}}$ with respect to 
$f_{m}$ as illustrated in 
Fig. \ref{slice_handcuff3} and a skein tree as illustrated in 
Fig. \ref{slice_handcuff3_2}. 
Then we have that 
\begin{eqnarray}\label{eq1}
a_{2}(K_{D}^{123}(f_{m}))&=&
a_{2}(J_{m})=a_{2}(J)-m-1 \\
&=&a_{2}(J_{m-1})+m+1-m-1 \nonumber\\
&=&a_{2}(J_{m-1})=\cdots=a_{2}(J_{0})=0. \nonumber
\end{eqnarray}
Then by Lemma \ref{vlemma2} and (\ref{eq1}) we have that 
\begin{eqnarray*}
v_{3}(K_{D}^{123}(f))&=&v_{3}(J_{m})=
v_{3}(J)+a_{2}(J_{m})+\frac{1}{2}(m+1)^{2}+1\\
&=&\left\{v_{3}(J_{m-1})-a_{2}(J_{m-1})-\frac{1}{2}(m+1)^{2}\right\}\\
&&+a_{2}(J_{m})+\frac{1}{2}(m+1)^{2}+1\\
&=&v_{3}(J_{m-1})+1=\cdots=v_{3}(J_{0})+m\\
&=&m.
\end{eqnarray*}
Since $K_{D}^{ij}(f_{m})$ is a trivial knot for 
any $1\le i<j\le 3$, 
we have that $\bar{n}_{123}(f_{m})=m$. 
Thus by Theorem \ref{main_inv2}, 
we have that 
$f_{m}$ is not delta edge-homotopic to the trivial spatial 
handcuff graph for any $m\in {\mathbb N}$, and 
$f_{i}$ and $f_{j}$ are 
not delta edge-homotopic for $i\neq j$. 
Note that the generalized version of $\bar{n}_{12}$ for $n=3$ 
as mentioned in 
Remark \ref{rem_inv} (2) vanishes for $f_{m}$ by (\ref{eq1}). 

\begin{figure}[htbp]
      \begin{center}
\scalebox{0.32}{\includegraphics*{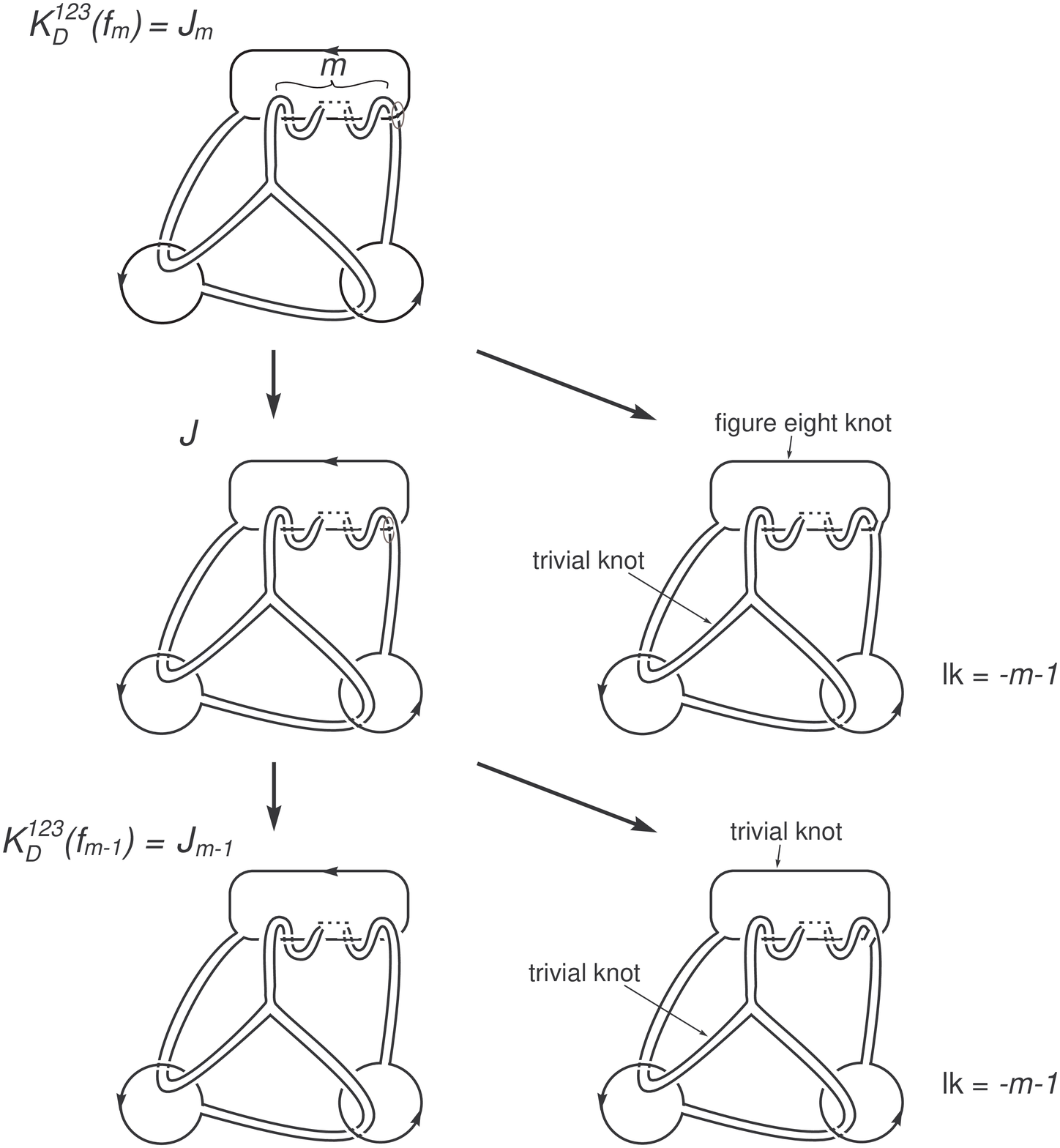}}
      \end{center}
   \caption{}
  \label{slice_handcuff3_2}
\end{figure} 
}
\end{Example}
\begin{Example}\label{ex4}
{\rm 
Let $f_{m}$ be the spatial handcuff graph for $m\in {\mathbb N}$ as illustrated 
in Fig. \ref{boundary_handcuff}. It is easy to see that $f_{m}$ is a 
$\partial$-spatial handcuff graph for any $m\in {\mathbb N}$. 
Since $L_{f_{m}}$ is a $2$-component boundary link, we have that 
$L_{f_{m}}$ is self $\Delta$-equivalent to the $2$-component trivial link 
for any $m\in {\mathbb N}$. 

\begin{figure}[htbp]
      \begin{center}
\scalebox{0.34}{\includegraphics*{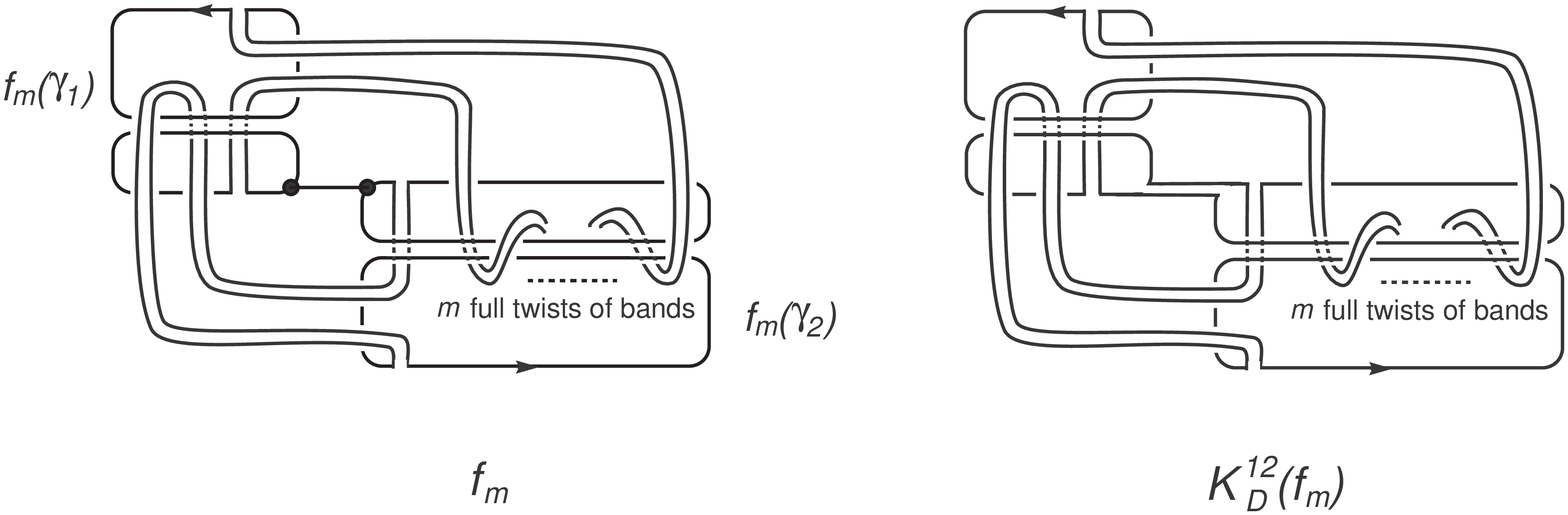}}
      \end{center}
   \caption{}
  \label{boundary_handcuff}
\end{figure} 

Now we consider the $D$-sum of $L_{f_{m}}$ with respect to $f_{m}$ 
as illustrated in 
Fig. \ref{boundary_handcuff}. 
Then by a calculation we have that $a_{2}(K_{D}^{12}(f_{m}))=-2m$. 
Since $f_{m}(\gamma_{1})$ and $f_{m}(\gamma_{2})$ are trivial knots, 
we have that $\bar{n}_{12}(f_{m})=-2m$. Thus by Theorem \ref{main_inv1}, 
we have that 
$f_{m}$ is not delta edge-homotopic to the trivial spatial 
handcuff graph for any $m\in {\mathbb N}$, and 
$f_{i}$ and $f_{j}$ are 
not delta edge-homotopic for $i\neq j$. 
}
\end{Example}

\begin{Example}\label{ex5}
{\rm 
Let $f$ be the spatial $3$-handcuff graph as illustrated 
in Fig. \ref{boundary_3_handcuff}. Note that 
$f|_{\gamma_{i}\cup\gamma_{j}\cup e_{i}\cup e_{j}}$ 
is the trivial spatial handcuff graph for any $1\le i<j\le 3$. 
It is easy to see that $f$ is a 
$\partial$-spatial $3$-handcuff graph. 
Since  $L_{f}$ is a $3$-component boundary link, 
we have that $\mu_{123}(L_{f_{m}})=0$ and $L_{f}$ is 
self $\Delta$-equivalent to the $3$-component trivial link. 
Note that $L_{f}$ is {\it Brunnian}, 
namely any $2$-component sublink of $L_{f}$ is trivial. 

\begin{figure}[htbp]
      \begin{center}
\scalebox{0.325}{\includegraphics*{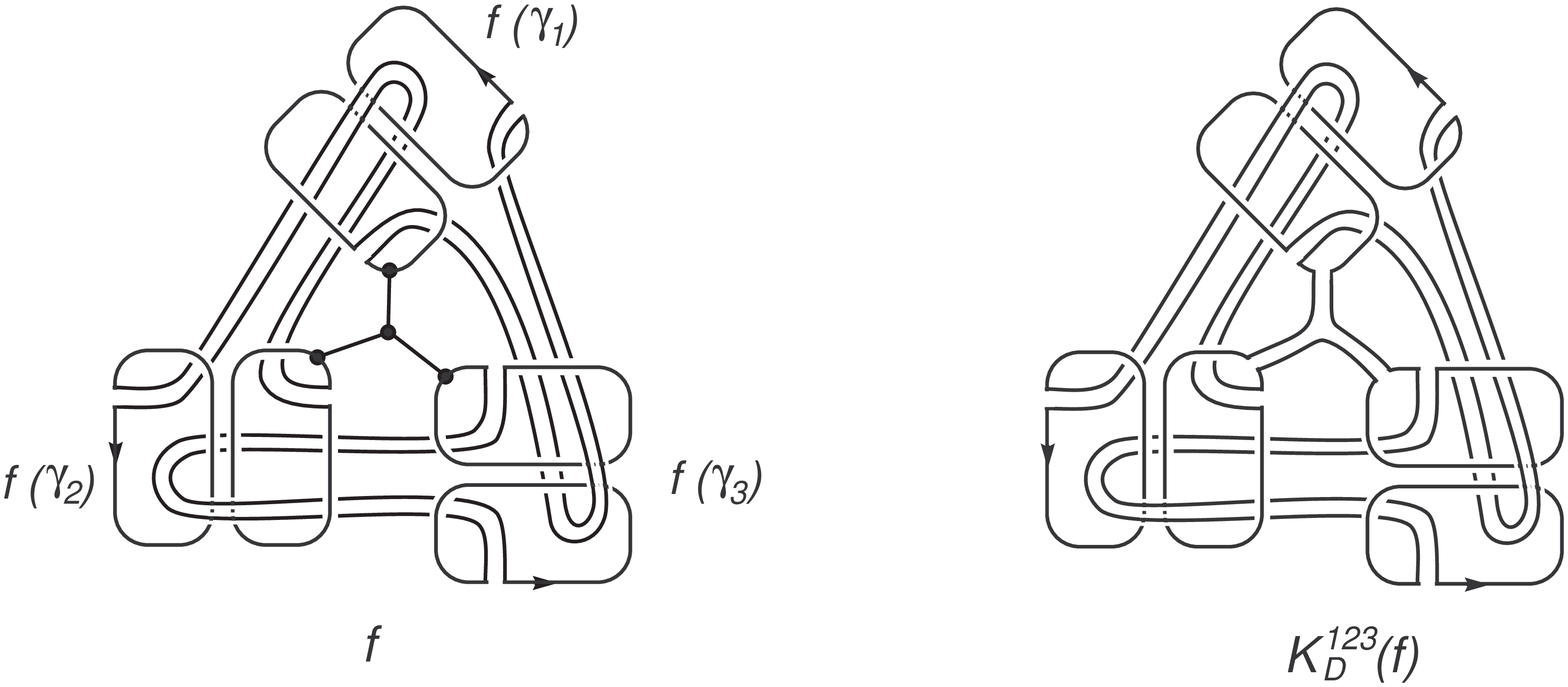}} 
      \end{center}
   \caption{}
  \label{boundary_3_handcuff}
\end{figure} 

Now we consider the $D$-sum of $L_{f}$ with respect to $f$ 
as illustrated in 
Fig. \ref{boundary_3_handcuff}. By a calculation we have 
that $a_{2}(K_{D}^{123}(f))=0$. 
Since $f(\gamma_{i})$ is a trivial knot for $i=1,2,3$, 
we have that 
the generalized version of $\bar{n}_{12}$ for $n=3$ 
as mentioned in 
Remark \ref{rem_inv} (2) vanishes for $f$. 
On the other hand, by a calculation we have that 
\begin{eqnarray*}
V_{K_{D}^{123}(f)}(t)&=&-t^{-12}+6t^{-11}-11t^{-10}+t^{-9}+28t^{-8}-52t^{-7}\\
&&+36t^{-6}+17t^{-5}-61t^{-4}+67t^{-3}-43t^{-2}+11t^{-1}\\
&&+22-57t+84t^{2}-78t^{3}+32t^{4}+23t^{5}-43t^{6}\\
&&+24t^{7}-4t^{9}-4t^{10}+5t^{11}+t^{12}-3t^{13}+t^{14}, \\
V_{K_{D}^{123}(f)}^{(3)}(1)&=&36. 
\end{eqnarray*}
Since $K_{D}^{ij}(f)$ is also a trivial knot for any $1\le i<j \le 3$, 
we have that $\bar{n}_{123}(f)=-1$. Thus by Theorem \ref{main_inv2},  
we have that $f$ is not delta edge-homotopic to the trivial spatial 
$3$-handcuff graph. 
}
\end{Example}

\vspace{0.2cm}
\begin{flushleft}
{\Large\bf Acknowledgments} 
\end{flushleft}

\vspace{0.2cm}
This work has been done while the author was visiting Waseda University and supported by Fellowship of the Japan Society for the Promotion of Science for Young scientists. The author is grateful to Professors Kouki Taniyama, 
Akira Yasuhara and Doctor Yukihiro Tsutsumi for their valuable comments. 
He would also like to thank Professor Kouji Kodama and his computer 
program KNOT \cite{kodamaknot} 
which enable the author to calculate the Conway polynomials and 
the Jones polynomials of knots which appear in section $4$ without 
difficulty.

{\normalsize
}

\end{document}